\documentclass[12pt, reqno]{amsart}
\usepackage{amsmath, amsthm, amscd, amsfonts, amssymb, graphicx, color}
\usepackage[bookmarksnumbered, colorlinks, plainpages]{hyperref}

\newtheorem{theorem}{Theorem}[section]
\newtheorem{lemma}[theorem]{Lemma}

\theoremstyle{definition}

\theoremstyle{remark}

\numberwithin{equation}{section}

\def\lsd{\mbox{$\mathop{\mathrel{\vrule height 5pt
 depth 0pt}\joinrel\mathrel\times} $}}

\title{Hypergroupoids and C*-algebras}
\author{Rohit Dilip Holkar}
\address{Universit\"at G\"ottingen, Mathematisches Institut 
Bunsenstr. 3-5 
D-37073 G\"ottingen, Germany }
\email{rohit@uni-math.gwdg.de}

\author{Jean Renault}
\address{Universit\'e d'Orl\'eans et CNRS  (UMR 7349 et FR2964), D\'epartement de Math\'ematiques, 
F-45067 Orl\'eans Cedex 2, France}
\email{Jean.Renault@univ-orleans.fr}

\thanks{The first author was supported by an Erasmus Mundus EuroIndia scholarship.}  


\usepackage{pdfsync}
\begin{document}
\maketitle
\markboth{Rohit Holkar and Jean Renault}
{Hypergroupoids and C*-algebras}

\begin{abstract}
Let $G$ be a locally compact groupoid. If $X$ is a free and proper $G$-space, then $(X*X)/G$ is a groupoid equivalent to $G$. We consider the situation where $X$ is proper but no longer free. The formalism of groupoid C*-algebras and their representations is suitable to attach C*-algebras to this new object.
\end{abstract}

\section{Introduction} This note stems from the elementary observation that the C*-category of a groupoid $G$ defined in \cite{ren:representations} can be extended from principal $G$-spaces to proper $G$-spaces. When $X$ is a principal locally compact $G$-space with invariant $r$-system $\alpha$, one can construct the $*$-algebra $(\alpha,\alpha)_c$ and its C*-completion $(\alpha,\alpha)$; it is the C*-algebra of the locally compact groupoid $(X*X)/G$ equipped with the Haar system induced by $\alpha$. When $X$ is only proper, the same formulas define the $*$-algebra $(\alpha,\alpha)_c$ and its C*-completion $(\alpha,\alpha)$; however $(X*X)/G$ is no longer a groupoid but a hypergroupoid (whose abstract definition is given in \cite{ren:induced}). Objects like $(X*X)/G$ generalize both hypergroups (when the $G$-space $X$ is transitive) and groupoids (when $X$ is free). While convolution algebras of measures are commonly associated with hypergroups, our construction gives convolution algebras of functions and C*-algebras. It also covers the construction of C*-algebras from Hecke pairs as in \cite{bc:hecke,kro:hecke,kro:hecke1}. In fact, the observation that $(X*X)/G$ is no longer a groupoid when $X$ is not a free $G$-space but that its convolution algebra can still be defined appears in this context (see \cite{lln:phase,cm:motives}). There, it is usual to introduce the reduced norm while the existence of a maximal norm is problematic. Our framework provides natural maximal and reduced norms on the hypergroupoids we consider. 

\section{The C*-category of a groupoid}

We review the framework and the main results of \cite{ren:representations}, but assuming that the $G$-spaces are proper and no longer free.  For the sake of simplicity, we consider here an untwisted groupoid $G$. Given a topological groupoid $G$ (with unit space $G^{(0)}$ and range and source maps $r$ and $s$), a left $G$-space is a topological space $X$ endowed with a continuous map $r_X: X\rightarrow G^{(0)}$, assumed to be open and onto, and a continuous action map $G*X\rightarrow X$, where $G*X$ is the subspace of composable pairs, i.e. $(\gamma,x)\in G\times X$ such that $s(\gamma)=r_X(x)$, sending $(\gamma,x)$ to $\gamma x$ in such a way that $(\gamma\gamma')x=\gamma(\gamma'x)$ for all composable triples $(\gamma,\gamma',x)$ and $ux=x$ if $u$ is a unit. One says that the $G$-space $X$ is proper [resp. free] if the map $G*X\rightarrow X\times X$ sending $(\gamma,x)$ to $(\gamma x,x)$ is proper [resp. injective]. We shall assume that the groupoid $G$ and the $G$-space $X$ are locally compact (and Hausdorff for the sake of simplicity).  If $X$ is a proper locally compact $G$-space, then the quotient space $X/G$ is locally compact. The image of $x\in X$ in the quotient space is denoted by $[x]$. If $X, Y$ are $G$-spaces, we endow $X*Y$ the space $X*Y$ of pairs $(x,y)$ such that $r_X(x)=r_Y(y)$ with the diagonal action $\gamma(x,y)=(\gamma x,\gamma y)$. It is  proper as soon as one of the factors is proper.  Given a $G$-space $X$, an $r_X$-system of measures is a family $\alpha=(\alpha^u)_{u\in G^{(0}}$ where $\alpha^u$ is a Radon measure on $X^u=r^{-1}(\{u\})$ with full support. We say that $\alpha$ is continuous if for all $f\in C_c(X)$, i.e. $f:X\rightarrow {\bf C}$ continuous with compact support, the function $u\mapsto \int f d\alpha^u$ is continuous on $G^{(0)}$. We say that $\alpha$ is invariant if for all $\gamma\in G$, $\gamma\alpha^{s(\gamma)}=\alpha^{r(\gamma)}$. The objects of our category are measured proper $G$-spaces $(X,\alpha)$, i.e. proper $G$-spaces $X$ endowed with a continuous invariant $r_X$-system of measures $\alpha$; we shall often omit the $G$-space $X$ and write $\alpha$ instead of $(X,\alpha)$. Before defining the C*-category $C^*(G)$, we first define the $*$-category (in the sense of \cite{glr:category}) $C_c(G)$: given two measured proper $G$-spaces $(X,\alpha)$ and $(Y,\beta)$, its set of arrows $(\alpha,\beta)_c$ consists of triples $(\alpha, f, \beta)$, where $f\in C_c((X*Y)/G)$. It is a complex vector space. Moreover, given measured proper $G$-spaces $(X,\alpha), (Y, \beta), (Z,\gamma)$ and $f\in C_c((X*Y)/G), g\in C_c((Y*Z)/G)$, we define
$$(\alpha, f, \beta)(\beta, g, \gamma)=(\alpha, f*_\beta g,\gamma)$$
where the convolution product is given by:
$$\leqno(1)\hskip 3cm f*_\beta g [x,z]=\int f[x,y]g[y,z] d\beta^{r_X(x)}(y)
$$
In this formula, a representative $(x,z)$ has been fixed and $[x,z]$ denotes its class. The integration is over a compact set because the map $\varphi^x: Y^{r_X(x)}\rightarrow (X*Y)/G$ defined by $\varphi^x(y)=[x,y]$ is proper. The resulting integral depends on $[x,z]$ only because of the invariance of $\beta$. One defines also
$$(\alpha, f, \beta)^*=(\beta, f^*, \alpha)$$
where the involution is given by 
$f^*[y,x]=\overline{f[x,y]}$.

\begin{lemma}(cf.\cite[Lemme 3.1]{ren:representations}) These operations are well defined and turn $C_c(G)$ into a $*$-category.
\end{lemma}

The next step is to define a C*-norm on the $*$-category $C_c(G)$. A unitary representation of $G$ is a pair $(m, H)$ where $m$ is a transverse measure class (\cite[Definition A.1.19]{adr:amenable}) and $H$ is a Borel $G$-Hilbert bundle. We recall that $m$ associates to $(X,\alpha)$ a measure class $m(\alpha)$ on $X/G$ in a coherent fashion. A unitary representation of $G$ defines by integration a representation of $C_c(G)$, that is, a functor into the W*-category of Hilbert spaces. It associates to the object $(X,\alpha)$ the Hilbert space $H(\alpha)=L^2(X/G,m(\alpha),X*H/G)$ and to the arrow $(\alpha,f,\beta)$ the operator $L(\alpha, f,\beta):H(\beta)\rightarrow H(\alpha)$ defined by
$$\langle \xi\sqrt\mu, L(\alpha, f,\beta)\eta\sqrt\nu\rangle=\int f[x,y]\langle\xi[x],\eta[y]\rangle\sqrt{(\mu\circ\dot\beta_1)(\nu\circ\dot\alpha_2)}[x,y]$$
where the sections $\xi\sqrt\mu\in H(\alpha)$ and $\eta\sqrt\nu\in H(\beta)$ are written as half-densities: $\mu$ [resp. $\nu$] is a measure on $X/G$ [resp. $Y/G$] in $m(\alpha)$ [resp. $m(\beta)$]. The systems of measures $\dot\beta_1$ and $\dot\alpha_2$ are induced by $\beta$ and $\alpha$ respectively as in \cite{ren:representations} or \cite[Lemma A.1.3]{adr:amenable} for the proper case. For example, one has $\int f d\dot\beta_1^{[x]}=\int f[x,y]d\beta^{r_X(x)}(y)$. By definition, the measures $m_1=\mu\circ\dot\beta_1$ and $m_2=\nu\circ\dot\alpha_2$ are equivalent; their geometric mean is the measure $(dm_1/dm_2)^{1/2}dm_2$. Note that by Cauchy-Schwarz inequality,
$$\|L(\alpha, f, \beta)\|\le \max\big(\sup_x\int |f[x,y]| d\beta^{r_X(x)}(y), \sup_y\int |f[x,y]| d\alpha^{r_Y(y)}(x)\big)$$
The I-norm of $f$ is defined as the left hand side. Just as in \cite{ren:representations}, we have

\begin{theorem}(cf. \cite[Proposition 3.5, Theorem 4.1]{ren:representations})
\begin{enumerate} 
\item Let $(m,H)$ be a unitary representation of a locally compact groupoid $G$. Then the above formulas define a representation $L$ of the $*$-category $C_c(G)$, called the integrated representation, which is continuous for the inductive limit topology and bounded for the I-norm.
\item Let $(G,\lambda)$ be a second countable locally compact groupoid with Haar system. Every  representation of the $*$-algebra $C_c(G,\lambda)$ in a separable Hilbert space which is non-degenerate and continuous for the inductive limit topology is equivalent to an integrated representation.
\end{enumerate}
\end{theorem}

We deduce from this theorem that, given a locally compact groupoid with Haar system $(G,\lambda)$, the $*$-category $C_c(G)$ can be completed into a C*-category by defining the full C*-norm $\|(\alpha,f,\beta)\|$ as the supremum of $\|L(\alpha,f,\beta)\|$ over all unitary representations of $G$ in separable Hilbert bundles. In particular, if $(X,\alpha)$ is a measured proper $G$-space, this defines the C*-algebra $(\alpha,\alpha)$. If moreover $X$ is a free $G$-space, $(X*X)/G$ is a groupoid equivalent to $G$; the  algebra $(\alpha,\alpha)$ is the full C*-algebra of this groupoid (endowed with the Haar system induced by $\alpha$) and is Morita equivalent to $C^*(G,\lambda)=(\lambda,\lambda)$. If $X$ is not free, $H=(X*X)/G$ is a hypergroupoid as defined in \cite{ren:induced}. Note that the full norm we have defined on the $*$-algebra $C_c(H)$ depends on $G$ since we consider only its representations which extend to the $*$-category $C_c(G)$. It is still true that $(\lambda,\alpha)$ is a full C*-module over $(\alpha,\alpha)$ but its algebra of compact operators is only an ideal of $C^*(G,\lambda)$. One has similar results with the regular representation and the reduced norm. If we identify $(G*X)/G=X$ through the map $(\gamma,x)\mapsto \gamma^{-1}x$, we obtain the various incarnations (2) of the formula (1): with $h\in C_c(G)$, $\xi,\eta\in C_c(X)$ and $f,g\in C_c((X*X)/G)$:
$$\leqno(2)\quad\begin{array}{cl} \xi f(y)=\int\xi(x)f[x,y]d\alpha(x)\quad&\langle\xi,\eta\rangle[x,y]=\int\overline{\xi(\gamma^{-1}x)}\eta(\gamma^{-1}y)d\lambda(\gamma)\\
h\xi(x)=\int h(\gamma)\xi(\gamma^{-1}x)d\lambda(\gamma)\quad&\langle\langle\xi,\eta\rangle\rangle(\gamma)=\int\xi(x)\overline{\eta(\gamma^{-1}x)}d\alpha(x)\\
f^*[y,x]=\overline{f[x,y]}\quad &f*g[x,z]=\int f[x,y]g[y,z]d\alpha(y) 
\end{array}$$

\section{Examples}

1. Let $K$ be a compact subgroup of a locally compact group $G$. The homogeneous space $X=G/K$ is a proper $G$-space equipped with an invariant measure $\alpha$. Then, $(X\times X)/G$ is the double coset hypergroup $H=K\backslash G/K$. The full and the regular representations of $G$ yield two C*-norms on the $*$-algebra of this hypergroup.  The resulting C*-algebras should be compared to the hypergroup C*-algebras which have been considered previously. The full C*-algebra $C^*(H)$ is defined in \cite{her:1992, her:1995} as the enveloping C*-algebra of the convolution algebra $L^1(H)$  (the authors are indebted to S. Echterhoff for drawing their attention to the work of P. Herman). There are examples of double coset hypergroups which have representations which do not extend to $G$ (this is carefully studied in \cite{her:1995}). Then the full norm coming from $G$ is strictly smaller than the full norm of the hypergroup. On the other hand, the reduced norms coincide since they come from the same regular representation.

2. Let $\Gamma_0$ be an almost normal subgroup of a discrete group $\Gamma$ as in \cite{bc:hecke,kro:hecke}. We equip $\Gamma/\Gamma_0$ with the counting measure. Since $\Gamma_0$ acts on $\Gamma/\Gamma_0$ with finite orbits, the convolution product is well-defined on $C_c(\Gamma_0\backslash \Gamma/\Gamma_0)$ which becomes the Hecke algebra ${\mathcal H}(\Gamma,\Gamma_0)$. Let $(G,K)$ be the Schlichting completion of $(\Gamma,\Gamma_0)$. Then ${\mathcal H}(\Gamma,\Gamma_0)$ can be identified with $C_c(K\backslash G/K)$ and we are in the situation of the first example.

3. A particular case of the next example, which generalizes the first example, is given in \cite[Section 1]{lln:phase}. Let $(G, \lambda)$ be a locally compact groupoid with Haar system and $K$ a closed subgroupoid with $K^{(0)}=G^{(0)}$. Assume that $K$ is a proper groupoid and that the map $r:G/K\rightarrow G^{(0)}$ has a $G$-invariant system of measures $\alpha$. Then $(X=G/K,\alpha)$ is a measured proper $G$-space. Thus we can construct the hypergroupoid $(X*X)/G=K\backslash G/K$ and its full and its reduced C*-algebras. If $K$ is principal, $(X*X)/G$ is a groupoid equivalent to $G$. The situation considered in \cite{lln:phase} is the case of a semi-direct groupoid $G=\Gamma\,\lsd Y$ where a group $\Gamma$ acts on a space $Y$ and $H=\Lambda\,\lsd Y$ where $\Lambda$ is a subgroup of $\Gamma$ acting properly on $Y$. The convolution algebra $C_c(K\backslash G/K)$ also appears in \cite{tk:biinvariant} (with $K$ compact), where the authors give a groupoid version of a Gelfand pair.




\section*{Acknowledgements}
The authors thank the GdR 2947 NCG which made this project possible.

\end{document}